\documentclass[12pt]{article}
\usepackage{amssymb}
\usepackage{latexsym}
\usepackage{amsmath}

\title{On $(n, k)$-extendable graphs and induced subgraphs} 

\author{Guizhen Liu\thanks{work supported  by NNSF of China and NSF of Shandong Province}\\ Department  of Mathematics, Shandong 
University,\\
Jinan, Shandong, P. R. China\\ \\Qinglin Yu\thanks{work supported by Natural Sciences and Engineering Research Council of Canada}\\ Center for Combinatorics, LPMC \\ Nankai University, Tianjin, China \\ and \\Department of Mathematics and Statistics
\\ Thompson Rivers University, Kamloops, BC, Canada}
\date{}
\begin{document}
\maketitle
\par
\begin{abstract}
Let $G$ be a graph with vertex set $V(G)$. Let $n$ and $k$ 
be non-negative integers such that $n + 2k \leq |V(G)| - 2$ 
and $|V(G)| - n$ is even. If when deleting any $n$ 
vertices of $G$ the remaining subgraph contains a matching of $k$ 
edges and every $k$-matching can be extended to a 1-factor, then $G$ is called an $(n, k)-${\bf extendable graph}. In this paper we present several results about $(n, k)$-extendable graphs and its subgraphs. In particular, we proved that if $G - V(e)$ is $(n, k)$-extendable graph for each $e \in F$ (where $F$ is a fixed 1-factor in $G$), then $G$ is $(n, k)$-extendable graph.

\end{abstract}
\par


{\bf AMS(1991) Subject Classification:} 05C70 \\

	Let $G$ be a simple graph with the vertex set $V(G)$ and the edge set 
$E(G)$. A {\bf matching} $M$ of $G$ is a subset of $E(G)$ such that any two edges of $M$ have no vertices in common. A matching of size $k$ is called a $k$-{\bf matching}. If $M$ is a matching so that every vertex (or except one) of $G$ is incident with an edge of $M$, then $M$ is called 1-{\bf factor} (or {\bf near} 1-{\bf factor}). 
  
	Let $S$ be a subset of $V(G)$. Denote by $G[S]$ the induced 
subgraph of $G$ by $S$ and we write $G - S$ for $G[V(G) \setminus S]$. $E(S, T)$ denotes the edges between two vertex sets $S$ and $T$.  The number of odd components of $G$ is denoted by $o(G)$.  

	Let $M$ be a matching of $G$. If there is a matching $M'$ of $G$ such that $M \subseteq M'$, then we say that 
$M$ can be extended to $M'$ or $M'$ is an extension of $M$. If each $k$-matching can be extended to a 1-factor, then $G$ is called $k$-{\bf extendable}.  A graph $G$ is called $n$-{\bf factor-critical} if 
after deleting any $n$ vertices the remaining subgraph of $G$ has a 1-factor. The properties of 2-factor-critical and $k$-extendable graphs were studied extensively by Lov\'asz and Plummer.  The history and applications of these topics can be found in \cite{key 5} and \cite{key 11}. Liu and Yu \cite{key 4} have introduced new concept, $(n, k)$-extendable graph, to combine the $n$-factor-criticality and the $k$-extendability.   

	Let $n$ and $k$ be non-negative integers 
such that $n + 2k \leq |V(G)| - 2$ and $|V(G)| - n$ is even. If
when deleting any $n$ vertices from $G$ the remaining subgraph of $G$ contains a $k$-matching and each $k$-matching in the subgraph can be extended to 1-factor, then $G$ is called a  $(n, k)$-{\bf extendable graph}.
Clearly, a graph is $(0, 0)$-extendable if and only if it has a 1-factor. 
Similarly, $(0, k)$-extendable graphs are precisely those $k$-extendable graphs and $(n, 0)$-extendable graphs are exactly $n$-critical graphs. A characterization and basic properties of $(n, k)$-extendable graphs were discussed in \cite{key 4}. 

	Nishimura and Saito \cite{key 7} and Yu \cite{key 13} studied the relationships between $k$-extendable graphs and its subgraphs and proved the followings\\

\noindent{\bf Theorem A.} (Nishimura and Saito \cite{key 7})  Let $G$ be a graph with a 1-factor.  If $G - V(e)$ is $k$-extendable for each $e \in E(G)$, then $G$ is $k$-extendable.\\

\noindent{\bf Theorem B.} (Yu \cite{key 13}) A graph $G$ is k-extendable if and only if for any matching M of size i
 ($1 \leq i \leq k$) the graph $G - V(M)$ is (k-i)-extendable.\\

	Based on Theorem B, Theorem A can be improved to the following:\\

\noindent{\bf Theorem 1.} Let $G$ be a graph with a 1-factor.  If $G - V(e)$ is k-extendable for each $e \in E(G)$ and $|V(G)| \geq 2k + 4$, then $G$ is ($k+1$)-extendable.

\noindent{\bf Proof:} Let $i = 1$ in Theorem B, then the result follows.  \hfill $\Box$\\

	In fact, the reverse of Theorem 1 is also true from Theorem B. Next we generalize this result to $(n, k)$-extendable graphs.\\

\noindent{\bf Theorem 2.} If $G - V(e)$ is an $(n, k)$-extendable graph for each $e \in E(G)$, then $G$ is $(n, k+1)$-extendable graph but may not be an $(n, k+2)$-extendable or $(n+2, k)$-extendable graph.

\noindent{\bf Proof:} Consider any vertex set $S$ and $(k+1)$-matching $M$ with $|S| = n$ and $V(M) \cap S = \emptyset$. Let $e$ be an edge of $M$.  Since $G - V(e)$ is $(n, k)$-extendable, there exists a 1-factor in $(G - V(e)) - (S \cup V(M - \{e\}) = G - (S \cup V(M))$.  Therefore, $G$ is an $(n, k+1)$-extendable graph. 

	To see that $G$ may not be $(n, k+2)$-extendable, we consider the graph
$$ H_1 = (2K_{2n+1}) + (K_n \cup (k+2)K_2)$$

	Then $H_1$ is not an $(n, k + 2)$-extendable graph by considering $S = V(K_n)$ and $(k+2)$-matching $(k+2)K_2$.  In the mean time, it is not hard to verify that for any $e \in E(H_1)$ $H_1 - V(e)$ is an $(n, k)$-extendable graph.

	Similarly, to see that $G$ may not be $(n+2, k)$-extendable, we consider the graph
$$ H_2 = (2K_{2n+1}) + (K_{n+2} \cup k K_2)$$

	Then $H_2$ is not an $(n+2, k)$-extendable graph but for any $e \in E(H_2)$ $H_2 - V(e)$ is an $(n, k)$-extendable graph.      
      \hfill $\Box$\\

	Before proceeding further, we quote two results from \cite{key 4} as lemmas.\\

\noindent{\bf Lemma 1.}  Let $G$ be an $(n, k)$-extendable graph.  Then it is also a $(n-2, k+1)$-extendable graph. \\

\noindent{\bf Lemma 2.} If G is an $(n, k)$-graph, then 

	(1) $G$ is also $(n -2 , k)$-extendable for $n \geq 2$;

	(2) $G$ is also $(n, k - 1)$-extendable for $k \geq 1$.\\

	For the convenience of the future arguments, we introduce one more term.  Let $S$ be a vertex set and $M$ a $k$-matching with $S \cap V(M) = \emptyset$.  If $G - S - V(M)$ has a 1-factor, then we say that G has a $(S, M)$-{\bf extension}.   

	Since an $(n+2, k)$-extendable or an $(n, k+2)$-extendable graph must be $(n, k+1)$-extendable, Theorem 2 indicates that $(n, k+1)$-extendability is the best possible under the general conditions.  But by introducing an additional condition on the size of graph in Theorem 2, we can improve it to the following: \\

\noindent{\bf Theorem 3.}  If $G - V(e)$ is an $(n, k)$-extendable graph ($n > 1$) for each $e \in E(G)$ and $V(G) \leq 2k + 3n + 4$, then $G$ is an $(n+2, k)$-extendable graph.

\noindent{\bf Proof:} Suppose that $G$ is not an $(n+2, k)$-extendable graph. By the definition, there exists a vertex set $S$ with $|S| = n+2$ and $k$-matching $M$ so that $G - S - V(M)$ has no 1-factor.  

	Let $G' = G - S - V(M)$.  From Tutte's Theorem, there exists a vertex set $S' \subseteq V(G')$ such that $o(G' - S') \geq |S'| + 2$.  \\

	{\it Claim 1.}  $G'- S'$ has exactly $|S'| + 2$ odd components.

	Otherwise, if $o(G' - S') \ne |S'| + 2$, by parity, then we have $o(G' - S') \geq |S'| + 4$.  Set $S_1 = S - \{a,b\}$ (where $a, b$ are any two vertices of $S$) and $S_1' = S' \cup \{a, b\}$.  Then 
$$o(G - S_1 - V(M) - S_1') = o(G - S - V(M) - S') = o(G' - S') \geq |S'| + 4 = |S_1'| +2 $$

	That is, $G - S_1 - V(M)$ has no 1-factor or $G$ has no $(S_1, M)$-extension.  But $|S_1| = n$  and $|M|= k$, so it contradicts to that G is $(n, k)$-extendable. \\

	{\it Claim 2.}  $S$ and $S'$ are independent sets. 

	If $S$ is not independent, let $e$ be an edge of $G[S]$ and $S_1 = S - V(e)$, then $G - V(e)$ has no $(S_1, M)$-extension.  This contradicts to the fact that $G - V(e)$ is an $(n, k)$-extendable graph. 
	
	Similarly, if $S'$ is not independent, let $e$ be an edge of $G[S']$, $S_1 = S - \{a, b\}$ (where $a$, $b$ are any two vertices of $S$) and $S_1' = S' - V(e) \cup \{a, b\}$, then $o(G - V(e)- S_1 - V(M) - S_1') = o(G - V(e) - S - V(M) - S') = o(G' - S') \geq |S_1'| + 2 $ or $G - V(e)$ has no $(S_1, M)$-extension.  This contradicts to that $G - V(e)$ is an $(n, k)$-extendable graph. \\	

	{\it Claim 3.}  $E(S, S') = \emptyset$.

	Otherwise, let $e = xy \in E(S, S')$ and $x \in S$, $y \in S'$.  Replacing the vertex $y$ by a vertex of $S - \{x\}$ and moving $y$ to $S$, then the new pair still have all of the properties of the old pair $S$ and $S'$ have but the new pair is against Claim 2, a contradiction. \\

	{\it Claim 4.}  No vertex in an even component is adjacent to $S \cup S'$.

	If there is an edge $e = xy$ so that $x \in S'$ and $y$ is in an even component. Set $S_1' = S' \cup \{y\}$.  Then  $$o(G - S - V(M) - S_1') = o(G - S - V(M) - S') +1  = o(G' - S') +1 \geq |S'| + 2 + 1 = |S_1'| + 2$$  
But $e = xy \in S_1'$, a contradiction to Claim 2. 
	
	Similarly, if there is an edge $e = xy$ so that $x \in S$ and $y$ is in an even component. Set $S_1' = S' -  \cup \{y\}$.  Then  $$o(G - S - V(M) - S_1') = o(G - S - V(M) - S') +1  = o(G' - S') +1 \geq |S'| + 2 + 1 = |S_1'| + 2$$  
But $e = xy \in E(S,S_1')$, a contradiction to Claim 3. \\

	With the preparation above, we can proceed to the proof of the theorem now. 

	From Theorem 2, G is $(n, k+1)$-extendable.  Applying Lemma 1 repeatedly we see that $G$ is $(\epsilon, (k+1+\lfloor n/2 \rfloor))$-extendable, where $\epsilon = 0$ or 1.  When $k$-matching $M$ is extended to a 1-factor (or near 1-factor) then $S \cup S'$ has to match to the vertices of odd components $\cup O_i$.  As $o(G' - S') = |S'| + 2$ and $n \geq 2$, so at least one of $O_i$'s has at least 3 vertices.  Choose an edge $e_1$ from such an odd component, say $O_1$, now we can extend $(k+1)$-matching $M \cup \{e_1\}$ to a 1-factor (or near 1-factor).  Thus $S \cup S'$ has to match to the vertices of $\cup O_i - V(e_1)$ and there exists an edge in $\cup O_i - V(e_1)$.  If this process is repeated, we can find $\lfloor n/2 \rfloor + 1$ disjoint edges in $\cup O_i$, namely, $\{e_1, e_2, \cdots, e_l\}$ (where $l = \lfloor n/2 \rfloor + 1$).  Since G is $(\epsilon, k + l)$-extendable,  $M \cup \{e_1, e_2, \cdots, e_l\}$ can be extended to a 1-factor (or near 1-factor), and thus $S \cup S'$ has to match to some vertices of $\cup O_i - V(e_1) - V(e_2)- \cdots - V(e_l)$.  Therefore, we have 
$$|V(G)| \geq 2|S \cup S'| + 2k + 2(\lfloor n/2 \rfloor + 1)$$
$$ \geq 2(n+2) + 2k + (n - 1) +2 = 2n + 4 + 2k + n + 1= 3n + 2k + 5$$
which contradicts to the given condition.  Hence, $G$ is an $(n+2, k)$-extendable graph. \hfill $\Box$\\

	Recently, Nishimura improved Theorem A by reducing the conditions required in the theorem.  Instead of checking the $k$-extendability of $G - V(e)$ for every edge $e$ in $G$, now one needs only checking the $k$-extendability of $G - V(e)$ for the edges belonging to a 1-factor of $G$.\\

\noindent{\bf Theorem C.} (Nishimura \cite{key 8})  Let $G$ be a graph with 1-factors and let $F$ be an arbitrary 1-factor of $G$.  If $G - V(e)$ is $k$-extendable graph (or $n$-factor-critical) for each $e \in F$, then $G$ is $k$-extendable (or $n$-factor-critical) graph.\\

	We will generalize the above result to $(n, k)$-extendable graphs.\\

\noindent{\bf Theorem 4.}  Let $G$ be a graph with 1-factors and let $F$ be an arbitrary 1-factor of $G$.  If $G - V(e)$ is $(n, k)$-extendable graph for each $e \in F$, then $G$ is $(n, k)$-extendable graph.

\noindent{\bf Proof:} We may assume that $n > 0$ and $k > 0$.

	We proceed to prove the theorem by contradiction.  Suppose that there exists a 1-factor $F$ of $G$ such that $G - V(e)$ is $(n, k)$-extendable for any $e \in F$ but $G$ is not $(n, k)$-extendable.  Then there exists a $k$-matching $M$ and a vertex set $S$ of size $n$, where $V(M) \cap S$ = $\emptyset$, such that $G - V(M)- S$ has no 1-factor.  Let $G'$ = $G - V(M)- S$. 
	Applying Tutte's 1-Factor Theorem, there exists $S' \subseteq V(G')$ so that $o(G' - S') > |S'|$.  By the parity, $o(G' - S') \geq |S'|+ 2$. Our aim is to find an edge $e \in F$ so that $G - V(e)$ is not $(n, k)$-extendable and thus leads to a contradiction.\\

	At first, we show that 1-factor $F$ can only match vertices from $V(M)$ to rest by the next claim.\\ 

	{\it Claim 1.} For the given $F$, $S$ and $G'$, we have 

		(i)  $F \cap E[S]$ = $\emptyset$;

		(ii) $F \cap E(S')$ = $\emptyset$;		

		(iii) $F \cap E(S, S')$ = $\emptyset$;

	To see (i), if $e \in F \cap E(S)$, then $|S - V(e)|$ = $n$ - 2 and $G - V(e)$ is not $(n - 2 , k)$-extendable.  Thus, $G$ is not $(n, k)$-extendable, a contradiction.

	To see (ii), if $e \in F \cap E(S')$, then $G' - V(e)$ has no 1-factor or $G - V(e)$ is not $(n, k)$-extendable, a contradiction.

	To see (iii), if $e \in F \cap E(S, S')$, where $e = ab$ and $a \in S$, $b \in S'$, choosing a vertex $c$ from an odd component of $G'- S'$ and then $S - \{a\} \cup \{c\}$ and $M$ can not be extended to a 1-factor as $o(G' - V(e) - S') > |S'| + 2 - 1$.

	From (i) - (iii), it follows that a 1-factor $F$ is in $E(S \cup V(M), G')$ or $E(S, V(M))$ or $E(G[V(M)])$. \\

	{\it Claim 2.}  $G'$ has no even components.
	
	Otherwise, let $D$ be an even component and let $e = ab$ be an edge of $F$, where $a \in V(D)$.
	
	If $b \in S$, choose $c \in V(D) - \{a\}$, then $T = S - \{b\}$ and $M$ can not extended to a 1-factor in $G - \{a, b\}$ as  $o((G' - V(e) - T - V(M))- S') \geq |S'|$ + 2, a contradiction.

	If $b \in V(M)$, consider an alternating path of $M \cup F$ with end-vertex $a$.  If another end-vertex $c$ of this alternating path is in $S$.  Similarly to the previous case, let $T$ = $S - \{c\} \cup \{x\}$ (where $x \in V(D) - \{a\}$ and $M' = M - \{bc'\} \cup \{ab\}$.  Then $G - \{c, c'\}$ (where $cc' \in F$) has no $(T, M')$-extension, a contradiction. 

	If $c$ is in $S'$, it is similar.

	If $c$ is in a component (either odd or even), let $T = S$ and $M' = M - \{bc'\} \cup \{ab\}$, then $G - \{c, c'\}$ has no $(T, M')$-extension as $G' - \{a, c\} - S'$ has at least $|S'| + 2$ odd components. \\

	{\it Claim 3.} $S'$ = $\emptyset$.

	If $S' \not= \emptyset$, let $a \in S'$, then $a$ is matched to a vertex $b$ in the 1-factor $F$ and $b$ must be in $V(M)$.  Consider an alternating path of $M \cup F$, say $abb' \cdots dd'c$. 

	If $c \in S'$, let $T = S$ and $M' = M - \{bb', dd'\} \cup \{ab, b'd\}$, then $G - \{d', c\}$ has no $(T, M')$-extension as $G' - \{a, c\}$ has no 1-factor.

	If $c \in S$, let $T = S - \{c\} \cup \{x\}$ (where $x$ is a vertex of a component) and $M' = M - \{bb', dd'\} \cup \{ab, b'd\}$, then $G - \{d', c\}$ has no $(T, M')$-extension as $G' - \{a, c\} - (S' - \{a\})$ has $o(G' - S')$ - 1 odd components, a contradiction.

	If $c \in C$ (where $C$ is any component), using the same argument we can see that $G' - \{a, c\} - (S' - \{a\})$ loses at most one odd component and obtain a contradiction. \\

	{\it Claim 4.}  $o(G' - S')$ = $o(G')$ = 2.

	Suppose $o(G') > 2$ (i.e., $o(G') \geq$ 4).  If there exists an edge $e \in F$ and $e \in E(S, C_1)$, choose $c$ from an odd component $C_2$, let $T$ = $S - \{b\} \cup \{c\}$ and $M' = M$, then $o(G' - \{a, c\}) \geq $ 2 or $G - \{a, b\}$ has no $(T, M')$-extension, a contradiction.

	Otherwise, all vertices in $\cup C_i$ are matched into $V(M)$.  Consider the alternating paths of $F \cup M$, there exists such a path starting with $C_i$ and ending $C_j$.  Let $c_ix_1y_1x_2y_2 \cdots x_my_mc_j$ be the alternating path, where $c_i \in C_i$, $c_j \in C_j$ and $c_ix_1$, $y_1x_2$, $\cdots$, $y_mc_j \in F$, $x_1y_1$, $x_2y_2$,  $\cdots$,  $x_my_m \in M$. 

	Let $T = S$ and $M' = M - \{ x_1y_1, \cdots, x_my_m \} \cup \{ y_1x_2, \cdots, y_mc_j \}$.  Then $G - \{ c_i, x_1 \}$ has no $(T, M')$-extension as $o(G' - \{c_i, c_j \}) \geq $ 2, a contradiction.\\

	{\it Claim 5.}  $F \cap E(S, V(M)) = \emptyset$.

	Consider the alternating path $ab \cdots c$ of $F \cup M$ with end-vertex $a$.  If $c \in S$, let $T = S - \{a, c \}$ and $M' = M - \{ bb' \} \cup \{cb' \}$, then $G - \{ a, b \}$ does not have $(T, M')$-extension, that is $G - \{ a, b \}$ is not ($n -2, k$)-extendable, a contradiction.  If $c \in C_1$ (where $C_1$ is an odd component) and $|C_1| \geq$ 3, choose $d \in V(C_1) - \{ c \}$ and let $T = S - \{ a \} \cup \{ d \}$ and $M' = M - \{ bb' \} \cup \{ b'c \}$. Then $G - \{ a, b \}$ (where $ab \in F$) has no $(T, M')$-extension as $o(G' - \{ c, d\}) \geq 2$.

	If $c \in C_1$ but $|C_1| $ = 1, then we have $|C_2| \geq 3$ because $G'$ has only two odd components, no even component and $|G'| \geq 4$. Suppose that $F \cap E(S, C_2) \not= \emptyset$.  Let $e = gh \in F \cap E(S, C_2)$, where $g \in V(C_2)$ and $h \in S$.  Choose $y \in V(C_2) - \{ g \}$ and set $T = S - \{ h \} \cup \{ y \}$ and $M' = M$, then $G - \{ g, h \}$ has no $(T, M')$-extension as $o(G' - \{g, y \}) \geq 2$, a contradiction.

	So we may assume $F \cap E(S, C_2) = \emptyset$.  In this case, all vertices of $C_2$ are matched to $V(M)$ in $F$.  Considering $F \cup M$, there must be an alternating path with both end-vertices in $V(C_2)$ or an alternating path starting in $V(C_2)$ and ending in $S$.  In either case, it yields a contradiction. \\

	Now we are ready to conclude the proof. \\

	Since $|S| \geq 1$ and $F \cap E(S, V(M)) = \emptyset$, there exists an edge $e = ab \in F$ from $S$ to an odd component $C_1$ (where $a \in S$, $b \in V(C_1)$). If $|C_1| \geq 3$, let $c  \in V(C_1) - \{ c \}$ and set $T = S - \{a\} \cup \{c\}$ and $M' = M$, then $G - \{a, b \}$ has no $(T, M')$-extension, a contradiction.  If $|C_1| = 1$, then $|C_2| \geq 3$.  Without loss of generality, we assume $F \cap E(S, C_2) = \emptyset$.  Thus, all vertices of $V(C_2)$ are matched to $V(M)$ in $F$.  Considering $F \cup M$, there exists an alternating path $P$ with both of ends in $C_2$ or an alternating path $P$ from $C_2$ to $S$.  

	Let $P = cx_1y_1d$, where $cx_1, y_1d \in F$ and $x_1y_1 \in M$.  If $c, d \in V(C_2)$, let $T = S$ and $M' = M - \{x_1y_1\} \cup \{dy_1\}$, then $G - \{c, x_1 \}$ has no $(T, M')$-extension as $o(G' - \{c, d\}) \geq 2$. If $c \in V(C_2)$ and $d \in S$, let $T = S - \{d\} \cup \{g\}$ (where $g \in V(C_2) - \{e\}$) and $M' = M - \{x_1y_1\} \cup \{dy_1\}$, then $G - \{c, x_1\}$ has no $(T, M')$-extension, a contradiction. 

	The proof is completed.     \hfill $\Box$\\

\par

\end{document}